\documentclass[10pt]{article}
\usepackage{amssymb,amsmath,amsfonts,amsthm,enumerate}
\usepackage{subfigure}
\usepackage{graphicx}
\usepackage{epsfig}
\usepackage{float}
\usepackage[usenames]{color}

\def\no{\noindent} \def\p{\partial} \def\nb{\nonumber}
\def\Vh0{\stackrel{\circ}{V}_h} 
   
\def\Om{\Omega}  \def\om{\omega}

\def\l{\label}  \def\f{\frac}  
\def\D{\end{document}}   
   
\def\m{\mbox}   \def\Box{\sharp}

\newcommand{\lc}
{\mathrel{\raise2pt\hbox{${\mathop<\limits_{\raise1pt\hbox
{\mbox{$\sim$}}}}$}}}

\newcommand{\gc}
{\mathrel{\raise2pt\hbox{${\mathop>\limits_{\raise1pt\hbox{\mbox{$\sim$}}}}$}}}

\newcommand{\ec}
{\mathrel{\raise2pt\hbox{${\mathop=\limits_{\raise1pt\hbox{\mbox{$\sim$}}}}$}}}

\def\bb{\begin{equation}}  \def\ee{\end{equation}}

\def\beqn{\begin{eqnarray}}  \def\eqn{\end{eqnarray}}

\def\beqnx{\begin{eqnarray*}} \def\eqnx{\end{eqnarray*}}

\def\bn{\begin{enumerate}} \def\en{\end{enumerate}}

\def\bd{\begin{description}} \def\ed{\end{description}}
\def\bfg{\begin{figure}} \def\efg{\end{figure}}

\newtheorem{lemma}{Lemma}[section]

\newtheorem{theorem}{Theorem}[section]
\newtheorem{assumption}{Assumption}[section]
\newtheorem{remark}{Remark}[section]

 \def\x{{\bf x}}

\title{Local Lipschitz Stability
for Inverse Robin Problems\\ in Some Elliptic and
Parabolic Systems}
\author{
Daijun Jiang\footnote{School of Mathematics and Statistics $\&$ Hubei Key Laboratory of Mathematical
Sciences, Central  China Normal University, Wuhan, 430079, P.R.China. The
work of this author was financially supported by
National Natural Science Foundation of China
(Nos.  11401241 and 11571265) and NSFC-RGC(China-Hong Kong, No. 11661161017).
(jiangdaijun@mail.ccnu.edu.cn)}
\and Jun
Zou\footnote{Department of Mathematics, The Chinese University of
Hong Kong, Shatin, Hong Kong. {\tt (zou@math.cuhk.edu.hk}). The work
of this author was substantially supported by Hong Kong RGC grants
(Projects 405110 and 404611). }}

\begin{document}
\maketitle

\begin{abstract}
In this work, we shall study the nonlinear inverse problems of recovering the Robin coefficients in
second order elliptic and parabolic systems with Robin boundary conditions, and establish their local Lipschitz stabilities.
 We shall first
 show the local Lipschitz stability for the elliptic inverse Robin problem and give some remarks to demonstrate why
 we slso need to consider the Robin  condition but not Neumann or Dirichlet condition on the accessible boundary.
 Then the new arguments are  generalized to help
establish a novel local Lipschitz stability for parabolic inverse Robin problems and some counterexamples shall also
be given to show why we also consider the Robin  condition on the accessible boundary.
\end{abstract}


%
\medskip
\section{Introduction}\l{sec:intro}
\setcounter{equation}{0}
We are concerned in this work with the determination of a spatially dependent Robin coefficient
in both stationary elliptic and time-dependent parabolic systems from
 measurement data on a partial boundary.
This is a highly nonlinear inverse problem and arises in several applications of practical importance.
The Robin coefficient may characterize the thermal properties of conductive materials
on the interface or certain physical processes near the boundary, e.g.,
it represents the corrosion damage profile
in corrosion detection (\cite{fang04}\cite {ing97}), and indicates the
thermal property in quenching processes (\cite {osm89}).

For the description of the model problems that are considered in this work,
we let $\Om= B_{r_2} (0)\backslash B_{r_1} (0) \subset R^2$ be an open bounded and connected annular domain,
where $0<r_1<r_2$ and $B_{r} (0)$ denotes a circle centered on the origin with radius $r>0$.
The boundary $\p\Om$ consists of two disjointed parts $\p\Om=\Gamma_i\cup\Gamma_a$, where
 $\Gamma_i=\p B_{r_1} (0)$ and $\Gamma_a=\p B_{r_2} (0)$ are respectively  the part
of the boundary that is inaccessible and accessible to experimental
measurements. Then we shall consider the inverse Robin problems associated with
the following elliptic boundary value problem
\bb \left\{ \begin{array} {rclll}
-\triangle u &=&f(\x) & \m{in} &\Om\,, \\
\frac{\p u}{\p n}+\gamma(\x) u&=&g(\x) & \m{on} &\Gamma_i\,, \\
\frac{\p u}{\p n}+ u &=&h(\x) & \m{on} &\Gamma_a\,,
\end{array}
\right. \label{ch1}
\ee
and parabolic initial boundary value problem
\bb \left\{ \begin{array} {rclll}
\p_t u-|\x|^2\triangle u &=&0 & \m{in} &\Om\times [0,T]\,, \\
\frac{\p u}{\p n}+\gamma(\x) u&=&g(\x,t) & \m{on} &\Gamma_i\times [0,T]\,, \\
\frac{\p u}{\p n}+ u &=&h(\x,t) & \m{on} &\Gamma_a\times [0,T]\,,\\
u(\x,0) &=&0 & \m{in} &\Om\,.
\end{array}
\right. \label{ch1t}
\ee
Functions $f$, $g$ and $h$ are the source strength,
ambient temperature and heat flux respectively. The coefficients $\gamma(\x)$ in (\ref{ch1}) and (\ref{ch1t}) represent
the Robin coefficients,  which
 will be the focus of our interest and are assumed to stay in the following feasible
constraint set:
$$
K:=\Big\{\gamma\in L^2(\Gamma_i);
\,\,0<\underline\gamma\leq\gamma(\x)\leq\bar\gamma \,\,a.e.\,\, {\rm on}\,\, \Gamma_i\Big\},
$$
where $\underline\gamma$ and $\bar\gamma$ are two positive constants.
For convenience, we often write the solutions of the systems (\ref{ch1}) and (\ref{ch1t}) as $u(\gamma)$ to emphasize
their dependence on the Robin coefficient $\gamma$.

We are now ready to formulate the inverse problems of our interest in this work.

{\bf Elliptic Inverse Robin Problem}: recover the Robin coefficient $\gamma(\x)$ in (\ref{ch1}) on the inaccessible
part $\Gamma_i$  from the measurable data $z$ of $u$ on the accessible part $\Gamma_a$.

{\bf Parabolic inverse Robin problem}: recover the Robin coefficient $\gamma(\x)$ in (\ref{ch1t})
on the inaccessible part $\Gamma_i$ from the measurable data $z$ of $u$ on the accessible part $\Gamma_a$
over the whole time range $[0,T]$.

Thanks for the unique continuation theorem \cite{isa06}, the uniqueness of the elliptic
and parabolic inverse Robin problems has been well established \cite{cha99}\cite{jiang15}\cite{jin10}.
The stability of the inverse Robin problems has also been studied for several years, but most of the studies are
 global logarithmic type stability estimate, see \cite{cha04}\cite{chou15}\cite{chou08}. There exists only a few results
on the Lipschitz stability of inverse Robin problems associated with the following elliptic equation:
\bb \left\{ \begin{array} {rclll}
-\triangle u &=&0 & \m{in} &D\,, \\
\frac{\p u}{\p n}+\gamma(\x) u&=&0 & \m{on} &\Gamma_1\,, \\
\frac{\p u}{\p n} &=&g(\x) & \m{on} &\Gamma_2\,,
\end{array}
\right. \label{j1}
\ee
where $D$ is  an open bounded and connected  domain in $R^2$ and $\Gamma_1$ and $\Gamma_2$ are the
 inaccessible and accessible part of the boundary $\p D$ respectively.
In \cite{cha99}, $Chaabane$ and $Jaoua$ proved a monotone global Lipschitz stability estimate. Simply speaking,
for $i = 1, 2$, let $u(\gamma_i)$ denote the solution of the boundary value problem \eqref{j1} corresponding
to $\gamma=\gamma_i$. If $\gamma_1\leq \gamma_2$, then we have an estimate of the form:
$\|\gamma_1-\gamma_2\|_1\leq C\|(u(\gamma_1)-u(\gamma_2))|_{\Gamma_a}\|_2$,
where $\|\cdot\|_1$ and $\|\cdot\|_2$ are some appropriate norms, which are different from
the standard Sobolev $H^1$ and $H^2$ norms. In \cite{chou04}, $Choulli$
established a local Lipschitz stability estimate for an arbitrary smooth domain without
the monotony condition $\gamma_1\leq \gamma_2$. The essential technique in the proof of the
local Lipschitz stability is the construction of a mapping, which is proved to be a $C^1$-diffeomorphism in
a neighborhood of a fixed element $\gamma_0\in K$. Very recently, $Hu$ and $Yamamoto$ also established in \cite{hu16} a global
H$\ddot{o}$lder stability estimate for an elliptic inverse Robin problem from
a single Cauchy data on an accessible boundary. The main arguments rely on the Schwarz reflection principle
with the Robin boundary condition and a novel interior estimate derived from
the elliptic Carleman estimate.

However, to our best knowledge, there are still no results available for the local Lipschitz stability for
parabolic inverse Robin problems. This will be one of the main novelties and contributions of the present work.
Another main novelty is that we construct many counterexamples creatively to show why we also consider
the Robin condition but not Neumann or Dirichlet condition on the accessible boundary for both the elliptic
Robin inverse problem and parabolic Robin inverse problem.
We shall first  prove the local Lipschitz stability of the
proposed elliptic inverse Robin problem and
the new arguments are then generalized to help us establish a novel local Lipschitz stability for the
inverse Robin problem associated with the  parabolic system \eqref{ch1t}, with some tricky and delicate
detailed modifications due to the complication of time dependence. It is worth of mention that it is the first time in literature to
establish the local Lipschitz stabilities for the inverse Robin problems associated with time-dependent
parabolic equations.

The rest of the paper is organized as follows. In Section\,\ref{sec:uni},
the local Lipschitz stability estimate
for the elliptic inverse Robin problem will be verified and  some remarks  to show why we also consider
the Robin condition but not Neumann or Dirichlet condition on the accessible boundary are given. In Section\,\ref{sec:un},
we shall establish a newly local Lipschitz stability estimate for the parabolic inverse Robin problem and
also give some remarks to show why we also consider
the Robin condition  on the accessible boundary

Throughout this work, $C$ is often used for a generic constant. We shall  write the norms of the spaces
$H^s(\Omega)$, $L^2(\Omega)$,
$H^{s}(\Gamma)$ and $L^2(\Gamma)$ (for some $\Gamma\subset \p\Omega$)  respectively
as $\|\cdot\|_{s, \Omega}$, $\|\cdot\|_{\Omega}$, $\|\cdot\|_{s, \Gamma}$ and
$\|\cdot\|_{\Gamma}$.

\section{Local Lipschitz stability for elliptic inverse Robin problem}\label{sec:uni}
\setcounter{equation}{0}
\setcounter{figure}{0}
\setcounter{table}{0}
In this section,  we shall establish the local Lipschitz stability for the proposed elliptic inverse Robin
problem.
We first give a preliminary lemma for recalling the classical well-posedness of the forward solution $u$ to system
(\ref{ch1}).

\begin{lemma}\label{lem:well}
{\rm (see \cite{evans10})\cite{gil01}} Let $\Om$ be an open bounded and connected domain with
$C^\infty$ boundary $\p\Om$,  $\gamma(\x)\in K$, $f(\x)\in L^2(\Om)$,
 $g(\x)\in L^{2}(\Gamma_i)$ and $h(\x)\in L^{2}(\Gamma_a)$,
then there exists a unique solution $u\in H^{2}(\Om)$ to system (\ref{ch1}) and it satisfies
\beqn
&&\|u\|_{2,\Om}
\leq C(\|f\|_{\Om}+\|g\|_{\Gamma_i}+
\|h\|_{\Gamma_a}).
\label{ch2}
\eqn

\end{lemma}

%


Then we study the differentiability of the solution $u(\gamma)$ to system \eqref{ch1} and
give its Fr$\acute{e}$chet derivative.

\begin{lemma}\label{lem:differential}
The solution $u(\gamma)$ of system \eqref{ch1} is continuously Fr$\acute{e}$chet differentiable and
its derivative  $u'(\gamma)d$ with direction $d\in L^\infty(\Gamma_i)$  solves the following system:

\bb \left\{ \begin{array} {rccll}
-\triangle u'(\gamma)d &=&0 & {\rm in} &\Om\,, \\
\frac{\p (u'(\gamma)d)}{\p n}+\gamma \,(u'(\gamma)d)&=&-d \,u(\gamma) & {\rm on} &\Gamma_i\,, \\
\frac{\p (u'(\gamma)d)}{\p n}+ u'(\gamma)d &=&0 & {\rm on} &\Gamma_a\,.
\end{array}
\right. \l{ch3}
\ee
\end{lemma}
\no {\it Proof}. For any $\gamma\in K$ and $d\in L^\infty(\Gamma_i)$ such that $\gamma+d\in K$, let
$v\equiv u(\gamma+d)-u(\gamma)-u'(\gamma)d$, then we have
\bb \left\{ \begin{array} {rccll}
-\triangle v &=&0 & \m{in} &\Om\,, \\
\frac{\p v}{\p n}+\gamma \,v&=&-d (u(\gamma+d)-u(\gamma)) & \m{on} &\Gamma_i\,, \\
\frac{\p v}{\p n}+v &=&0 & \m{on} &\Gamma_a\,.
\end{array}
\right.
\ee
From estimate \eqref{ch2} and the Sobolev embedding theorem, we obtain
\beqnx
\|v\|_{1,\Om}&\leq& C\|d (u(\gamma+d)-u(\gamma))\|_{\Gamma_i}\leq
C\|d\|_{L^\infty(\Gamma_i)}\|u(\gamma+d)-u(\gamma)\|_{\f{1}{2},\Gamma_i}\\
&\leq& C\|d\|_{L^\infty(\Gamma_i)}\|u(\gamma+d)-u(\gamma)\|_{1,\Om}.
\eqnx
As $\psi\equiv u(\gamma+d)-u(\gamma)$ satisfies the following elliptic equation
\bb \left\{ \begin{array} {rccll}
-\triangle \psi &=&0 & \m{in} &\Om\,, \\
\frac{\p \psi}{\p n}+\gamma \,\psi&=&-d u(\gamma+d) & \m{on} &\Gamma_i\,, \\
\frac{\p \psi}{\p n}+\psi &=&0 & \m{on} &\Gamma_a\,.
\end{array}
\right.
\ee
Similarly, we can show from estimate \eqref{ch2} and the Sobolev embedding theorem that
\beqnx
\|\psi\|_{1,\Om}\leq C\|d u(\gamma+d)\|_{\Gamma_i}\leq
C\|d\|_{L^\infty(\Gamma_i)}\|u(\gamma+d)\|_{\f{1}{2},\Gamma_i}
\leq C\|d\|_{L^\infty(\Gamma_i)}\|u(\gamma+d)\|_{1,\Om}\leq C\|d\|_{L^\infty(\Gamma_i)}.
\eqnx
Thus it follows directly that
\beqnx
\f{\|u(\gamma+d)-u(\gamma)-u'(\gamma)d\|_{1,\Om}}{\|d\|_{L^\infty(\Gamma_i)}}
\rightarrow 0~~{\rm as}~~\|d\|_{L^\infty(\Gamma_i)}\rightarrow 0,
\eqnx
which means that $u(\gamma)$ is Fr$\acute{e}$chet differentiable and $u'(\gamma)d$ is its derivative.

Next, we verify the continuity of $u'(\gamma)d$. Let $\phi\in L^\infty(\Gamma_i)$, then
$y\equiv u'(\gamma+\phi)d-u'(\gamma)d$ satisfies
\bb \left\{ \begin{array} {rccll}
-\triangle y &=&0 & \m{in} &\Om\,, \\
\frac{\p y}{\p n}+\gamma \,y&=&-d(u(\gamma+\phi)-u(\gamma))-\phi u(\gamma+\phi) & \m{on} &\Gamma_i\,, \\
\frac{\p y}{\p n}+y &=&0 & \m{on} &\Gamma_a\,.
\end{array}
\right.
\ee
Once again estimate \eqref{ch2} implies that
\beqnx
\|y\|_{1,\Om}&\leq& C(\|d(u(\gamma+\phi)-u(\gamma))\|_{\Gamma_i}+\|\phi u(\gamma+\phi)\|_{\Gamma_i})\\
&\leq& C\|d\|_{L^\infty(\Gamma_i)}\|u(\gamma+\phi)-u(\gamma)\|_{1,\Om}+C\|\phi\|_{L^\infty(\Gamma_i)}\| u(\gamma+\phi)\|_{\Gamma_i}\\
&\leq&C\|d\|_{L^\infty(\Gamma_i)}\|\phi\|_{L^\infty(\Gamma_i)}+C\|\phi\|_{L^\infty(\Gamma_i)},
\eqnx
which tends to 0 when $\|\phi\|_{L^\infty(\Gamma_i)}$ tends to 0.
$\Box$

Now let $\gamma^*\in K$ be the true Robin coefficient for the proposed elliptic Robin inverse problem.
Then from Lemma \ref{lem:differential} we see that $u(\gamma^*)$  is continuously Fr$\acute{e}$chet differentiable and
its derivative  $w(d)\doteq u'(\gamma^*)d$ satisfies the following elliptic system
\bb \left\{ \begin{array} {rccll}
-\triangle w(d) &=&0 & \m{in} &\Om\,, \\
\frac{\p w(d)}{\p n}+\gamma^* \,w(d)&=&-d u(\gamma^*) & \m{on} &\Gamma_i\,, \\
\frac{\p w(d)}{\p n}+w(d) &=&0 & \m{on} &\Gamma_a\,.
\end{array}
\right. \label{ch4}
\ee
With the help of $\om(d)$, we define a bounded and linear operator from $L^2(\Gamma_i)$ to $L^2(\Gamma_a)$ as follows:
\beqnx
N(d)=\f{\p w(d)}{\p n}~~{\rm on}~~{\Gamma_a},~~\forall\,d\in L^2(\Gamma_i).
\eqnx
For establishing the local Lipschitz stability estimate, we need to make some assumptions.

\begin{assumption}\label{ass1}
$u(\gamma^*)\neq 0$ almost everywhere on $\Gamma_i$.
\end{assumption}
This assumption is very natural, as one can see from the second equation of system \eqref{ch1} that if
$u(\gamma^*)= 0$ a.e. on $\Gamma_i$ then the true Robin coefficient $\gamma^*$ is not identifiable.

\begin{assumption}\label{ass2}
When we consider the polar coordinates $(r,\,\theta)$, functions $f(\x),\,g(\x),\,h(\x)$ and true Robin
coefficient $\gamma^*$ are only dependent on $r$ but independent of $\theta$.
\end{assumption}

\begin{lemma}\label{lem:sur}
Under Assumptions \ref{ass1}-\ref{ass2}, the operator $N$ is bijective and $\|N^{-1}\|\leq C$.
\end{lemma}
\no {\it Proof}. We first show $N$ is injective, i.e., if $N(d)=0$ then $d=0$. Indeed, if
\beqnx
N(d)=\f{\p w(d)}{\p n}=0~~{\rm on}~~{\Gamma_a},
\eqnx
then the third equation of \eqref{ch4} implies that $w(d)=0$ on $\Gamma_a$. Therefore,
we know that $w(d)$ is also the solution of the following system
\bb \left\{ \begin{array} {rccll}
-\triangle w &=&0 & \m{in} &\Om\,, \\
w &=&0 & \m{on} &\Gamma_a\,, \\
\frac{\p w}{\p n} &=&0 & \m{on} &\Gamma_a\,,
\end{array}
\right.
\ee
which implies that $w(d)=0$ in $\Om$ by the unique continuation principle \cite{isa06}  and thus
$d u(\gamma^*)=0$ on $\Gamma_i$. Then we easily get $d=0$ on $\Gamma_i$ by Assumption \ref{ass1}.

Now we prove $N$ is surjective, which means that
for any element $\varphi=\frac{\p w}{\p n}\in L^2(\Gamma_a)$, we want to seek $d\in L^2(\Gamma_i)$ such
that $N(d)=\varphi$. To do so, we use the separation of variables in polar coordinates to solve \eqref{ch4} and
thus the first equation transforms to the following equation with variables $r$ and $\theta$
\bb
\frac{\p^2 w}{\p r^2}+\frac{1}{r}\frac{\p w}{\p r}+\frac{1}{r^2}\frac{\p^2 w}{\p \theta^2}=0.
\label{jd1}
\ee
The Assumption \ref{ass2} implies that the solutions $u$ and $w$ of system \eqref{ch1} and \eqref{ch4} respectively
are all just dependent on $r$ but independent on $\theta$ and thus we have
\bb
\frac{\p^2 w}{\p r^2}+\frac{1}{r}\frac{\p w}{\p r}=0.
\label{jd2}
\ee
It is easy to get the general solutions of \eqref{jd2} that
\bb
w=c_1+c_2\ln r,
\label{jd3}
\ee
where $c_1$ and $c_2$ are two real number. Hence we get on $\Gamma_a$, i.e., $r=r_2$ that
\beqnx
&&\varphi=\frac{\p w}{\p n}=\frac{\p w}{\p r}=c_2\frac{1} {r_2},\\
&&w=c_1+c_2\ln r_2,
\eqnx
and thus
\beqnx
0=\frac{\p w}{\p n}+w=c_2\frac{1}{r_2}+(c_1+c_2\ln r_2)=c_1+c_2(\frac{1}{r_2}+\ln r_2).
\eqnx
For simplicity, we let $c_2=-1$, then $c_1=\frac{1}{r_2}+\ln r_2$ and obtain
\beqnx
w=\frac{1}{r_2}+\ln r_2-\ln r ~~{\rm in}~~\Om.
\eqnx
Hence, we have on $\Gamma_i$, i.e., $r=r_1$ that
\beqnx
&&-d u(\gamma^*)(r_1)=\frac{\p w}{\p n}+\gamma^*(r_1) \,w=-\frac{\p w}{\p r}+\gamma^*(r_1) \,w\\
&=&\frac{1}{r_1}+\gamma^*(r_1)\,(\frac{1}{r_2}+\ln r_2-\ln r_1),
\eqnx
which with Assumption \ref{ass1} implies that
\beqnx
d=-\frac{1}{r_1u(\gamma^*)(r_1)}-\frac{\gamma^*(r_1)}{u(\gamma^*)(r_1)}\,(\frac{1}{r_2}+\ln r_2-\ln r_1).
\eqnx

As $N$ is linear, bounded and bijective, then by the Open Mapping Theorem \cite{con03}, we know that
 $N^{-1}$ exists and $\|N^{-1}\|$ is bounded, i.e., there exists a positive constant $C$ such that
 \beqnx
 \|N^{-1}\|\leq C.~~~~~\Box
 \eqnx
Finally,  for convenience, we shall write for any positive constant $b$ that
\beqnx
N(\gamma^*,b)=\{\gamma\in K;\,\,
\|\gamma-\gamma^*\|_{\Gamma_i}\leq b\}.
\eqnx
We are now ready to establish the local Lipschitz stability for elliptic inverse Robin problem.

\begin{theorem}\label{thm:locale}
Under Assumptions \ref{ass1}-\ref{ass2}, there exists a positive constant $b$ such that the
following stability estimate holds:
\bb
\|u(\gamma_1)-u(\gamma_2)\|_{\Gamma_a}\geq C\|\gamma_1-\gamma_2\|_{\Gamma_i},~~\forall~\gamma_1,\gamma_2\in N(\gamma^*,b).
\label{ch5}
\ee

\end{theorem}
\no {\it Proof}. We first introduce an important mapping
\beqnx
\theta:\gamma\in L^2(\Gamma_i)\rightarrow \f{\p u(\gamma)} {\p n}\in L^2(\Gamma_a),
\eqnx
which is continuously Fr$\acute{e}$chet-differentiable from Lemma \ref{lem:differential} and it is obviously to see

\beqnx
\theta'(\gamma^*)d=\f{\p u'(\gamma^*)d}{\p n}= N(d).
\eqnx
Then it follows from lemma \ref{lem:sur} that $\theta'(\gamma^*)^{-1}= N^{-1}$ and

\beqnx
\|\theta'(\gamma^*)^{-1}\|=\| N^{-1}\|
\leq C.
\eqnx
By  the inverse function theorem \cite{dep03} we find that $\theta(\gamma^*)$ is $C^1$-diffeomorphism
on a neighborhood of $\gamma^*$, consequently $\theta(\gamma^*)^{-1}$ is locally Lipschitz continuous and
 $\|(\theta(\gamma^*)^{-1})'\|=\|\theta'(\gamma^*)^{-1}\|\leq C$.
Thus, there exists a neighborhood  $N(\gamma^*,b)$ of $\gamma^*$ such that  the Lipschitz
constant is less equal to $2\|(\theta(\gamma^*)^{-1})'\|$, i.e., for any $\gamma_1,\gamma_2\in N(\gamma^*,b)$,
it holds that

\beqn
\|\gamma_1-\gamma_2\|_{\Gamma_i}&\leq& 2\|(\theta(\gamma^*)^{-1})'\|
\|\f{\p u(\gamma_1)} {\p n}-\f{\p u(\gamma_2)} {\p n}\|_{\Gamma_a}\nb\\
&\leq&2C\|u(\gamma_2) - u(\gamma_1) \|_{\Gamma_a},
\label{cch3}
\eqn
where we used the fact that $\f{\p u(\gamma)} {\p n}=-u(\gamma)$ on $\Gamma_a$ in the second inequality. $\Box$

\section{Local Lipschitz stability for parabolic inverse Robin problem}\label{sec:un}
\setcounter{equation}{0}
\setcounter{figure}{0}
\setcounter{table}{0}
In this section, we shall establish the
local Lipschitz stability for the proposed parabolic inverse Robin problem and  give some
counterexamples to show why we also consider the Robin condition on accessible part.
We first give a preliminary lemma for recalling the classical well-posedness of the forward solution $u$ to system
(\ref{ch1t}).

\begin{lemma}\label{lem:wellt}
{\rm (see \cite{evans10} \cite{jin12})} Let $\Om$ be an open bounded and connected domain with
$C^\infty$ boundary $\p\Om$,    $\gamma(\x)\in K$,
 $g(\x,t)\in L^2(0,T;L^{2}(\Gamma_i))$ and  $h(\x,t)\in L^2(0,T;L^{2}(\Gamma_a))$,
then there exists a unique solution $u\in L^2(0,T;H^{2}(\Om))$ to system (\ref{ch1t}) and it satisfies

\beqn
\|u\|_{L^2(0,T;H^{2}(\Om))}
\leq C(\|g\|_{L^2(0,T;L^{2}(\Gamma_i))}+
\|h\|_{L^2(0,T;L^{2}(\Gamma_a))}).
\label{c1}
\eqn

\end{lemma}

Now we study the differentiability of the solution $u(\gamma)$ to system \eqref{ch1t} and
give its Fr$\acute{e}$chet derivative.

\begin{lemma}\label{lem:differentialp}
The solution $u(\gamma)$ of system \eqref{ch1t} is continuously Fr$\acute{e}$chet differentiable and
its derivative $u'(\gamma)d$ with direction $d\in L^\infty(\Gamma_i)$  solves the following system:

\bb \left\{ \begin{array} {rclll}
\p_t (u'(\gamma)d)-|\x|^2\triangle (u'(\gamma)d) &=&0 & {\rm in} &\Om\times [0,T]\,, \\
\frac{\p (u'(\gamma)d)}{\p n}+\gamma \,(u'(\gamma)d)&=&-d\,u(\gamma) & {\rm on} &\Gamma_i\times [0,T]\,, \\
\frac{\p (u'(\gamma)d)}{\p n}+u'(\gamma)d &=&0 & {\rm on} &\Gamma_a\times [0,T]\,,\\
(u'(\gamma)d)(\x,0) &=&0 & {\rm in} &\Om\,,
\end{array}
\right. \label{c2}
\ee
\end{lemma}
\no {\it Proof}. For any $\gamma\in K$ and $d\in L^\infty(\Gamma_i)$ such that $\gamma+d\in K$, let
$v\equiv u(\gamma+d)-u(\gamma)-u'(\gamma)d$, then we have
\bb \left\{ \begin{array} {rclll}
\p_t v-|\x|^2\triangle v &=&0 & \m{in} &\Om\times [0,T]\,, \\
\frac{\p v}{\p n}+\gamma \,v&=&-d (u(\gamma+d)-u(\gamma)) & \m{on} &\Gamma_i\times [0,T]\,, \\
\frac{\p v}{\p n}+v &=&0 & \m{on} &\Gamma_a\times [0,T]\,,\\
v(\x,0) &=&0 & \m{in} &\Om\,,
\end{array}
\right.
\ee
From estimate \eqref{c1} and the Sobolev embedding theorem, we obtain
\beqnx
\|v\|_{L^2(0,T;H^1(\Om))}&\leq& C\|d (u(\gamma+d)-u(\gamma))\|_{L^2(0,T;L^2(\Gamma_i))}\\
&\leq&
C\|d\|_{L^\infty(\Gamma_i)}\|(u(\gamma+d)-u(\gamma))\|_{L^2(0,T;H^{\f{1}{2}}(\Gamma_i))}\\
&\leq& C\|d\|_{L^\infty(\Gamma_i)}\|(u(\gamma+d)-u(\gamma))\|_{L^2(0,T;H^1(\Om))}.
\eqnx
As $\psi\equiv u(\gamma+d)-u(\gamma)$ satisfies the following parabolic equation
\bb \left\{ \begin{array} {rclll}
\p_t \psi-|\x|^2\triangle \psi &=&0 & \m{in} &\Om\times [0,T]\,, \\
\frac{\p \psi}{\p n}+\gamma \,\psi&=&-d u(\gamma+d) & \m{on} &\Gamma_i\times [0,T]\,, \\
\frac{\p \psi}{\p n}\psi &=&0 & \m{on} &\Gamma_a\times [0,T]\,,\\
\psi(\x,0) &=&0 & \m{in} &\Om\,,
\end{array}
\right.
\ee
Similarly, we can show from estimate \eqref{c1} and the Sobolev embedding theorem that
\beqnx
\|\psi\|_{L^2(0,T;H^1(\Om))}&\leq& C\|d u(\gamma+d)\|_{L^2(0,T;L^2(\Gamma_i))}\leq
C\|d\|_{L^\infty(\Gamma_i)}\|u(\gamma+d)\|_{L^2(0,T;H^{\f{1}{2}}(\Gamma_i))}\\
&\leq& C\|d\|_{L^\infty(\Gamma_i)}\|u(\gamma+d)\|_{L^2(0,T;H^1(\Om))}\leq C\|d\|_{L^\infty(\gamma_i)}.
\eqnx
Thus it follows directly that
\beqnx
\f{\|u(\gamma+d)-u(\gamma)-u'(\gamma)d\|_{L^2(0,T;H^1(\Om))}}{\|d\|_{L^\infty(\Gamma_i)}}
\rightarrow 0~~{\rm as}~~\|d\|_{L^\infty(\Gamma_i)}\rightarrow 0,
\eqnx
which means that $u(\gamma)$ is Fr$\acute{e}$chet differentiable and $u'(\gamma)d$ is its derivative.

Next, we verify the continuity of $u'(\gamma)d$. Let $\phi\in L^\infty(\Gamma_i)$, then
$y\equiv u'(\gamma+\phi)d-u'(\gamma)d$ satisfies
\bb \left\{ \begin{array} {rclll}
\p_t y-|\x|^2\triangle y &=&0 & \m{in} &\Om\times [0,T]\,, \\
\frac{\p y}{\p n}+\gamma \,y&=&-d(u(\gamma+\phi)-u(\gamma))-\phi u(\gamma+\phi) & \m{on} &\Gamma_i\times [0,T]\,, \\
\frac{\p y}{\p n}+y &=&0 & \m{on} &\Gamma_a\times [0,T]\,,\\
y(\x,0) &=&0 & \m{in} &\Om\,,
\end{array}
\right.
\ee
Once again estimate \eqref{c1} implies that
\beqnx
\|y\|_{1,\Om}&\leq& C(\|d(u(\gamma+\phi)-u(\gamma))\|_{L^2(0,T;L^2(\Gamma_i))}+\|\phi u(\gamma+\phi)\|_{L^2(0,T;L^2(\Gamma_i))})\\
&\leq&C\|d\|_{L^\infty(\Gamma_i)}\|\phi\|_{L^\infty(\Gamma_i)}+C\|\phi\|_{L^\infty(\Gamma_i)},
\eqnx
which tends to 0 when $\|\phi\|_{L^\infty(\Gamma_i)}$ tends to 0.
$\Box$

Now let $\gamma^*\in K$ be the true Robin coefficient for the proposed parabolic Robin inverse problem,
then from Lemma \ref{lem:differentialp} we see that $u(\gamma^*)$  is continuously Fr$\acute{e}$chet differentiable and
its derivative  $w(p)\doteq u'(\gamma^*)p$ satisfies the following parabolic system
\bb \left\{ \begin{array} {rclll}
\p_t w(p)-|\x|^2\triangle w(p) &=&0 & \m{in} &\Om\times [0,T]\,, \\
\frac{\p w(p)}{\p n}+\gamma^* \,w(p)&=&-pu(\gamma^*) & \m{on} &\Gamma_i\times [0,T]\,, \\
\frac{\p w(p)}{\p n}+w(p) &=&0 & \m{on} &\Gamma_a\times [0,T]\,,\\
w(p)(\x,0) &=&0 & \m{in} &\Om\,.
\end{array}
\right. \l{c3}
\ee

Next we define a bounded and linear operator from $L^2(\Gamma_i)$ to $L^2(0,T;L^2(\Gamma_a))$ as follows:
\beqn
N(p)=\p w(p)/\p n~~{\rm on}~~\Gamma_a\times [0,T],~~\forall\,p\in L^2(\Gamma_i).
\l{c7}
\eqn
For establishing the local Lipschitz stability estimate, we need to make some assumptions.

\begin{assumption}\label{ass3}
$u(\gamma^*)\neq 0$ almost everywhere on $\Gamma_i\times [0,T]$.
\end{assumption}
This assumption is very natural, as one can see from the second equation of system \eqref{ch1t} that if
$u(\gamma^*)= 0$ a.e. on $\Gamma_i\times [0,T]$ then the true Robin coefficient $\gamma^*$ is not identifiable.

\begin{assumption}\label{ass4}
The solution $u(\gamma^*)$ to system \eqref{ch1t} is completely separated into its spatial and temporal components
and has the formulation $u(\gamma^*)=te^{t-\ln t}v(\x)$.
When we consider the polar coordinates $(r,\,\theta)$, $v(\x)$ and true Robin
coefficient $\gamma^*$ are only dependent on $r$ but independent of $\theta$.
\end{assumption}

\begin{lemma}\label{lem:sur2p}
Under Assumptions \ref{ass3}-\ref{ass4}, $N$ is bijective and $\|N^{-1}\|$ is bounded.
\end{lemma}
\no {\it Proof}. We first show $N$ is injective, i.e., if $N(p)=0$ then $p=0$. Indeed, if
\beqnx
N(p)=\f{\p w(p)}{\p n}=0~~{\rm on}~~{\Gamma_a}\times [0,T],
\eqnx
then the third equation of \eqref{c3} implies that $w(p)=0$ on $\Gamma_a\times [0,T]$. Therefore,
we know that $w(p)$ is also the solution of the following system
\bb \left\{ \begin{array} {rclll}
\p_t w(p)-|\x|^2\triangle w(p) &=&0 & \m{in} &\Om\times [0,T]\,, \\
w(p)&=&0 & \m{on} &\Gamma_a\times [0,T]\,, \\
\frac{\p w(p)}{\p n} &=&0 & \m{on} &\Gamma_a\times [0,T]\,,\\
w(p)(\x,0) &=&0 & \m{in} &\Om\,.
\end{array}
\right.
\ee
which implies that $w(p)=0$ in $\Om\times [0,T]$ by the unique continuation principle \cite{isa06}  and thus
$p u(\gamma^*)=0$ on $\Gamma_i\times [0,T]$. Then we easily get $p=0$ on $\Gamma_i$ by Assumption \ref{ass3}.

Now we prove $N$ is surjective, which means that
for any element $\varphi=\frac{\p w}{\p n}\in L^2(0,T;L^2(\Gamma_a))$, we want to seek $p\in L^2(\Gamma_i)$ such
that $N(p)=\varphi$. To do so, we first see from Assumption \ref{ass4} that $w(p)$ should also have the expression
$w(p)=te^{t-\ln t}V(\x)$ and when we consider the polar coordinates $(r,\,\theta)$, $V(\x)$ is only dependent on $r$.

Next, we use the separation of variables in polar coordinates to solve \eqref{c3} and
thus get
\beqnx
0=\p_t w(p)-|\x|^2\triangle w(p)&=&te^{t-\ln t}
(V-r^2\frac{\p^2 V}{\p r^2}-r\frac{\p V}{\p r}-\frac{\p^2 V}{\p \theta^2})\nb\\
&=&te^{t-\ln t}
(V-r^2\frac{\p^2 V}{\p r^2}-r\frac{\p V}{\p r}).
\eqnx
This gives that
\beqnx
r^2\frac{\p^2 V}{\p r^2}+r\frac{\p V}{\p r}-V=0,
\eqnx
which is an Euler ordinary differential equation and its general solutions, for any constants $c_1,c_2$,
\beqnx
V=c_1r+c_2r^{-1}.
\eqnx

 Hence we obtain on $\Gamma_a\times [0,T]$, i.e., $r=r_2$ that
\beqnx
&&\varphi=\frac{\p w}{\p n}=\frac{\p w}{\p r}=te^{t-\ln t}(c_1-c_2\frac{1}{r_2^2}),\\
&&w=te^{t-\ln t}(c_1r_2+c_2\frac{1}{r_2}),
\eqnx
and thus
\beqnx
0=\frac{\p w}{\p n}+w=te^{t-\ln t}\{c_1(1+r_2)-c_2(\frac{1}{r_2^2}-\frac{1}{r_2})\}.
\eqnx
For simplicity, we let $c_2=1$, then $c_1=\frac{1-r_2}{r_2^2(1+r_2)}$ and obtain
\beqnx
w=te^{t-\ln t}(\frac{1-r_2}{r_2^2(1+r_2)}r+\frac{1}{r}) ~~{\rm in}~~\Om.
\eqnx
Therefore, with Assumption \ref{ass4}, we have on $\Gamma_i\times[0,T]$, i.e., $r=r_1$ that
\beqnx
&&-d\, te^{t-\ln t}v(r_1)=\frac{\p w}{\p n}+\gamma^*(r_1) \,w=-\frac{\p w}{\p r}+\gamma^*(r_1) \,w\\
&=&te^{t-\ln t}\{-\frac{1-r_2}{r_2^2(1+r_2)}+\frac{1}{r_1^2}+\gamma^*(r_1)(\frac{1-r_2}{r_2^2(1+r_2)}r_1+\frac{1}{r_1})\},
\eqnx
which with Assumption \ref{ass3} implies that
\beqnx
d=-\frac{1}{v(r_1)}\{\frac{r_2-1}{r_2^2(1+r_2)}+\frac{1}{r_1^2}+\gamma^*(r_1)(\frac{1-r_2}{r_2^2(1+r_2)}r_1+\frac{1}{r_1})\}.
\eqnx

As $N$ is linear, bounded and bijective, then by the Open Mapping Theorem \cite{con03}, we know that
 $N^{-1}$ exists and $\|N^{-1}\|$ is bounded, i.e., there exists a positive constant $C$ such that
 \beqnx
 \|N^{-1}\|\leq C.~~~~~\Box
 \eqnx

We are now ready to establish the local Lipschitz stability for parabolic inverse Robin problems.
\begin{theorem}\label{thm:locale}
Under Assumption \ref{ass3}-\ref{ass4},  we have
\bb
\|\gamma_1-\gamma_2\|_{\Gamma_i}\leq C\|u(\gamma_1)-u(\gamma_2)\|_{L^2(0,T;L^2(\Gamma_a))}.
\label{c6}
\ee
\end{theorem}
\no {\it Proof}.  We first introduce an important mapping
\beqnx
\theta:\gamma\in L^2(\Gamma_i)\rightarrow \f{\p u(\gamma)} {\p n}\in L^2(0,T;L^2(\Gamma_a)),
\eqnx
which is continuously Fr$\acute{e}$chet-differentiable from Lemma \ref{lem:differentialp} and get

\beqnx
\theta'(\gamma^*)p=\f{\p u'(\gamma^*)p}{\p n}= N(p).
\eqnx
Then it follows from lemma \ref{lem:sur2p} that $\theta'(\gamma^*)^{-1}= N^{-1}$ and

\beqnx
\|\theta'(\gamma^*)^{-1}\|=\| N^{-1}\|
\leq C.
\eqnx
By  the inverse function theorem \cite{dep03} we find that $\theta(\gamma^*)$ is $C^1$-diffeomorphism
on a neighborhood of $\gamma^*$, consequently $\theta(\gamma^*)^{-1}$ is locally Lipschitz continuous and
 $\|(\theta(\gamma^*)^{-1})'\|=\|\theta'(\gamma^*)^{-1}\|\leq C$.
Thus, there exists a neighborhood  $N(\gamma^*,b)$ of $\gamma^*$ such that  the Lipschitz
constant is less equal to $2\|(\theta(\gamma^*)^{-1})'\|$, i.e., for any $\gamma_1,\gamma_2\in N(\gamma^*,b)$,
it holds that

\beqn
\|\gamma_1-\gamma_2\|_{\Gamma_i}&\leq& 2\|(\theta(\gamma^*)^{-1})'\|
\|\f{\p u(\gamma_1)} {\p n}-\f{\p u(\gamma_2)} {\p n}\|_{L^2(0,T;L^2(\Gamma_a))}\nb\\
&\leq&2C\|u(\gamma_2) - u(\gamma_1) \|_{L^2(0,T;L^2(\Gamma_a))}.~~~~~~~\Box
\eqn

\begin{remark} \label{re:4}
Assumption \ref{ass4} can be improved to more general situation:  the solution $u(\gamma^*)$ to system \eqref{ch1t} is
completely separated into its spatial and temporal components
and has the formulation $u(\gamma^*)=F(t)v(\x)$ with $F(0)=0$, $F(t)\neq 0$ and $F'(t)=F(t)$.
When we consider the polar coordinates $(r,\,\theta)$, $v(\x)$ and true Robin
coefficient $\gamma^*$ are only dependent on $r$ but independent of $\theta$.
\end{remark}

\end{document}